\documentclass[a4paper, 11pt, reqno]{amsart}
\usepackage[utf8]{inputenc}
\usepackage[T1]{fontenc}
\usepackage[frenchb,english]{babel}
\usepackage{enumerate}
\usepackage{dsfont, amsthm, amsfonts, amsmath, amssymb, stmaryrd, latexsym}
\usepackage{amscd}
\usepackage{array}
\usepackage{color}
\usepackage{graphics}
\usepackage{multicol}
\usepackage{longtable}
\usepackage{graphics}
\usepackage{tikz}
\usepackage{xypic}
\usetikzlibrary{arrows}
\usepackage{capt-of}
\usepackage{geometry}
\usepackage{float}
\usepackage{url}		
\theoremstyle{plain}
\newtheorem{thm}{Theorem}
\newtheorem{prop}{Proposition}
\newtheorem{cor}{Corollary}
\newtheorem{lem}{Lemma}

\theoremstyle{remark}
\newtheorem{rem}{Remark}

\theoremstyle{definition}
\newtheorem{Def}{Definition}

\def\N{\mathds{N}}
\def\R{\mathds{R}}
\def\Q{\mathbb{Q}}
\def\C{\mathbb{C}}
\def\Z{\mathbb{Z}}

\begin{document}
\selectlanguage{english}
\title[On  quaternion algebras that split...]{On  quaternion algebras that split over  specific  quadratic number fields}

\author[V. Acciaro]{Vincenzo Acciaro}
\author[D. Savin]{Diana Savin}
\author[M. Taous]{Mohammed Taous}
\author[A. Zekhnini]{Abdelkader Zekhnini}

\address{Vincenzo Acciaro, Dipartimento di Economia, Universit\`a di Chieti--Pescara,
Viale della Pineta, 4,  65127 Pescara, Italy}
\email{v.acciaro@unich.it}

\address{Diana Savin, Faculty of Mathematics and Computer Science, Ovidius University,
Bd. Mamaia 124, 900527, Constanta, Romania}
\email{savin.diana@univ-ovidius.ro; dianet72@yahoo.com}

\address{Mohammed Taous, Department of  Mathematics, Faculty of  Sciences and Technology, Moulay Ismail University of Meknes, Errachidia, Morocco.}
\email{taousm@hotmail.com}

\address{Abdelkader Zekhnini, Mohammed Premier University, Pluridisciplinary Faculty, Department Mathematics and Informatics, Nador, Morocco}
\email{zekha1@yahoo.fr}

\begin{abstract}
Let $d$ and $m$ be two distinct  squarefree  integers and $\mathcal{O}_K$ the ring of integers of the  quadratic field $K=\mathbb{Q}(\sqrt{d})$. Denote by  $ H_K(\alpha, m)$  a quaternion algebra  over $K$, where
 $\alpha\in \mathcal{O}_K$. In this paper we give necessary and sufficient conditions for  $ H_K(\alpha, m)$ to split over $K$  for some values of  $\alpha$,  and we obtain a complete characterization of division quaternion algebras $ H_K(\alpha, m)$  over
 $K$ whenever $\alpha$ and  $m$ are two distinct positive prime integers. Examples are given  involving prime Fibonacci numbers.
\end{abstract}
\maketitle
\section{Introduction}
Let $K$ be a field of characteristic $\neq 2$. If $\alpha, \beta\in K^*=K\backslash\{0\}$, then there exists  a unique unital associative
$K$-algebra of dimension $4$ with $K$-basis  $\left\{ 1, i, j, k\right\}$ such that  $i^2 = \alpha$, $j^2 = \beta$ and $ij = -ji=k$.
 This $K$-algebra will be denoted  by  $ H_K(\alpha, \beta)$. Its presentation, as a $K$-algebra, is  given by $K\{i, j\}/(i^2-\alpha, j^2-\beta, ij=-ji)$. A quaternion algebra over $K$ is a $K$-algebra isomorphic to such an algebra for
some $\alpha, \beta\in K^*$. The classical instance, where $K=\R$, is  $ H_K(-1, -1)=\mathbb{H}$, the  Hamilton's quaternions $(\alpha=\beta=-1)$.

The classification of quaternion algebras over $K$ can be rephrased in terms of
quadratic forms, and a more detailed description depends on the field $K$. In this
vein, the most important  question one  may ask about a quaternion algebra $ H_K(\alpha, \beta)$ is whether it is
isomorphic to the matrix ring $\mathrm{M}_2(K)$;  if so, we say that  $ H_K(\alpha, \beta)$  splits over $K$. For
example, every quaternion algebra over $\C$ (or an algebraically closed field)  splits,
and a quaternion algebra $ H_{\R}(\alpha, \beta)$ over $\R$  splits if and only if $\alpha > 0$ or $\beta> 0$ (see \cite{Mi}).

Quaternion algebras  are central simple
algebras over $K$ (i.e.  associative and non-commutative algebras without two sided ideals whose center is precisely
$K$, see \cite{Mi}) of dimension $4$ over $K$. Recall that the dimension $d$ of a central simple algebra
$A$ over a field $K$ is always a perfect square, and its square root $n$ is defined to be
the degree of $A$.

The theory of central simple algebras (and thus in particular quaternion algebras and
cyclic algebras) has strong connections with algebraic number theory, combinatorics, algebraic geometry, coding theory, computer science and signal theory (see \cite{Gi-Sz, Ja, Mi}). Quaternion algebras have been  studied in many papers that deal with conditions for an algebra to be split or  a division
algebra (e.g. \cite{Sav16, Sav17, Sav14}).

To decide whether a quaternion algebra is a division algebra or it splits, different approaches  are used: quadratic forms,  the associated
conics (which are  projective plane curves defined by the homogeneous equation $\alpha x^2 + \beta y^2 = z^2$), cyclotomic fields, $p$-adic fields and some other properties of associative algebras (e.g. \cite{Sav16, Sav17, Sav14}).

In this paper we adopt two distinct approaches. The first consists of studying  the ramification of
certain integral primes, and we obtain a nice characterization of quaternion division
algebras $H_K(p, q)$ solely in terms of quadratic residues, assuming that $p$ and $q$ are
positive primes and $K$ is  quadratic number field. The second one, makes use of the Hilbert symbol  in order to
    obtain necessary and sufficient conditions for a quaternion algebra
    $ H_K(\alpha, m)$ to split over a quadratic field $K=\Q(\sqrt d)$,  for some  integers $\alpha$  of $\mathcal{O}_K$, the ring of integers of the  quadratic field $K$, where $m$  and $d$ are  squarefree  integers.
\section{Preliminaries}
Let us collect  some results  that we will used in the squeal. We begin by recalling the definition of the ring of integers of a number field.
\begin{Def} Let $K$ be a number field.
	\begin{enumerate}[\indent\rm1.]
		\item The ring of integers      $\mathcal{O}_{K}$  of  $K$ is defined as follows
$$\mathcal{O}_{K}=\left\{x\in K:\;P(\alpha)=0\text{ for some monic polynomial } P\in\Z[X]\right\};$$
		\item $\mathcal{O}_{K}^\times$, the unit group of $K$, is the set of  invertible elements of $\mathcal{O}_{K}$;
		\item If $K=\Q(\sqrt{d})$ is a quadratic field, where $d>0$ is a  squarefree  integer, then there exists  a  unit $\varepsilon_d$, called the fundamental unit of $K$, and it is a generator (modulo the roots of unity) of the unit group  $\mathcal{O}_{K}^\times$.
	\end{enumerate}
\end{Def}

 Let $K$ be a number field and  $q$  an element of the  quaternion algebra $ H_K(\alpha, \beta)$ over $K$; then we can write $q = a + bi + cj + ek$, where $a$, $b$, $c$ and $e$ are in $K$.
The conjugate $\overline{q}$ of $q$  is defined as $\overline{q} = a -bi - cj - ek$.    The norm map is defined by $N$: $q\mapsto N(q) = q\overline{q}=a^2-\alpha b^2-\beta c^2+\alpha\beta e^2\in K$.
So  $N$ may be regarded as a quadratic form in the four variables
$a, b, c$ and  $e$.
\begin{Def}
A non-zero vector $q\in H_K(\alpha, \beta)$ is said to be isotropic if $N(q) =0$. 
\end{Def}
\begin{Def}\label{1}
An associative algebra $A$ over a field is called a division algebra if and only if it has a multiplicative identity element $1 \neq 0$ and every non-zero element  in $A$ has a left and a right multiplicative inverses.	If $A$ is a finite-dimensional algebra, then $A$ is a division algebra if and only if $A$  has
no nontrivial zero divisors.
\end{Def}
\begin{Def}
A central simple algebra $A$ over a field $K$ is called split by $K$ if $A$ is isomorphic to the full matrix algebra $M_2(K)$ over $K$.
\end{Def}
\begin{thm}[\cite{Gi-Sz, Vig80, Voi}]\label{twodoteightone} Let   $K$   be a field with   char $K \neq 2$   and let  $\alpha, \beta \in K \backslash \{  0 \}$. Then  the  quaternion algebra
$H= H_K(\alpha, \beta)$ is either  split  or  a division algebra. Furthermore,  the following statements are equivalent:
\begin{enumerate}[\indent\rm1.]
	\item $H \cong  H_K(1, 1)\cong M_2(K)$;
	\item $H$ is not a division algebra;
	\item $H$ has an isotropic element;
    \item the equation $\alpha x^2 + \beta y^2 = 1$ has a solution $(x, y)\in K\times K$;
    \item $\alpha\ \ ($resp. $\beta)$ is a norm from $K(\sqrt{\beta})\ \  ($resp. from $K(\sqrt{\alpha}))$.
\end{enumerate}
\end{thm}
The stamements above can be checked by using    local techniques. For a number field $K$, it is well known that its ring of integers  $\mathcal{O}_K$ is a Dedekind ring. If   $\alpha\in K$ and       $\mathcal{P}$    is a  prime ideal  of $\mathcal{O}_K$, we have
$\alpha\mathcal{O}_K=\mathcal{P}^{v_{\mathcal{P}}(\alpha)}I,$ where $v_{\mathcal{P}}(\alpha)$ is the highest power of $\mathcal{P}$  dividing $\alpha\mathcal{O}_K$ and $I$ is an ideal of $\mathcal{O}_K$.  The map $\alpha\mapsto N(\mathcal{P})^{-v_{\mathcal{P}}(\alpha)}$ is the non-Archimedean absolute value,
where $N$ denotes the absolute norm map.
 We denote by $K_{\mathcal{P}} $, the completion of $K$ with respect to this absolute value.
We conclude this section by defining the concept of ramification of a quaternion algebra  at a prime ideal.
\begin{Def}\label{ghent3}
Let  $\mathcal{P}$ be a prime ideal of $\mathcal{O}_K$. We say that ${H}_K(\alpha, \beta)$ is ramified at $\mathcal{P}$ if ${H}_{K_\mathcal{P}}(\alpha, \beta)$ is a division ring. The set of ramified primes of ${H}_K(\alpha, \beta)$ will be denoted by $\mathrm{Ram}({H}_K(\alpha, \beta))$. The reduced discriminant $D_{H_K(\alpha, \beta)}$ of the quaternion algebra  ${H}_K(\alpha, \beta)$
is defined as the product of those prime ideals of $\mathcal{O}_{K}$ that
 ramify in  ${H}_K(\alpha, \beta)$.
\end{Def}

\section{Symbols over number fields}
If $K$ is a number field of degree $n$, with signature $(r, s$) and $L$  is a quadratic extension of $K$, then $L=K(\sqrt{\alpha})$ for some  $\alpha\in \mathcal{O}_K$. If $\mathcal{P}$ is a prime ideal of $\mathcal{O}_K$, then
the behavior of its  extension $\mathcal{P}\mathcal{O}_L$ to $L$  is as follows:
$$
\mathcal{P}\mathcal{O}_L=\begin{cases}\mathcal{B}_1\mathcal{B}_2,&\mathcal{P} \mbox{ splits into $2$ different prime ideals of $\mathcal{O}_L$};\\
        \mathcal{P},&\mbox{$\mathcal{P}$ remains prime ($\mathcal{P}$ is called inert in $L$)};\\
   \mathcal{B}^2, & \mbox{$\mathcal{P}$     ramifies in $L$}.
 \end{cases}
$$
Note that $\mathcal{B}_i\cap K=\mathcal{B}\cap K=\mathcal{P}$ and $\mathcal{P}\cap \Q=p\Z$, where $p$ is a prime number. We say that $\mathcal{B}$ (resp. $\mathcal{B}_i$)  lies  above $\mathcal{P}$. The following properties hold:

\begin{itemize}
	\item $\mathcal{O}_K/\mathcal{P}$ is a finite field and $N(\mathcal{P}):=\mid\mathcal{O}_K/\mathcal{P}\mid=p^{f_{\mathcal{P}}}$ with $f_{\mathcal{P}}\in\N$. $N(\mathcal{P})$ is called the absolute norm of $\mathcal{P}$ and $f_{\mathcal{P}}$ is its inertia degree.
			\item If $v_{\mathcal{P}}(\alpha)$ is odd, then $\mathcal{P}$ is ramified in $L=K(\sqrt \alpha)$.
			\item If $v_{\mathcal{P}}(\alpha)$ is even and $\mathcal{P}$ does not divide $2$, then $\mathcal{P}$ is unramified in $L=K(\sqrt \alpha)$.
			\item The infinite primes are just the $r+s$ non equivalent archimedean absolute values coming  from the $r$ real and $s$ pairs of complex emebeddings.
		\end{itemize}
{\bf I. The quadratic residue symbol}.\\
The quadratic residue symbol   $\left(\dfrac{\alpha}{\mathcal{P}}\right)$  is classically defined as follows. Let  $\mathcal{P}$ be a prime ideal of
$\mathcal{O}_K$.  If $\alpha$ is a square in $K$, we let  $\left(\dfrac{\alpha}{\mathcal{P}}\right)=1$. Otherwise, we let
$$\displaystyle\left(\dfrac{\alpha}{\mathcal{P}}\right)=\left\{
\begin{array}{rl}
1, & \hbox{if $\mathcal{P}$ splits in $K(\sqrt \alpha)$;} \\
-1, & \hbox{if $\mathcal{P}$ is inert in $K(\sqrt \alpha)$;} \\
0, & \hbox{if $\mathcal{P}$ ramifies in $K(\sqrt \alpha)$.}
\end{array}
\right.
$$
\begin{thm}[\cite{Gr03, Lm00}]
	If $\mathcal{P}$ is prime to  the ideal generated by $\alpha$ and to the ideal generated by  $2$, then $$\displaystyle\left(\dfrac{\alpha}{\mathcal{P}}\right)\equiv \alpha^{\frac{N(\mathcal{P})-1}{2}}\pmod{\mathcal{P}}.$$
\end{thm}
\begin{rem}\label{ghent1}~\
		\begin{enumerate}[\indent\rm1.]
\item  If $K=\Q$, then $\mathcal{P}=p\Z$ where $p$ is  a prime integer. In this case the symbol $\left(\dfrac{\alpha}{\mathcal{P}}\right)$ is denoted by $\left(\dfrac{\alpha}{p}\right)$ and called the Legendre symbol.
If $\left(\dfrac{\alpha}{p}\right)=1$ (in this case, we say that $\alpha$ is a quadratic residue modulo $p$) and $p\equiv 1\pmod 4$, then
$\left(\dfrac{\alpha}{p}\right)_4$ will denote the rational biquadratic symbol which is equal to
$1$ or $-1$, according as $(\alpha)^{\frac{p-1}{4}}\equiv 1 \text{ or } -1\pmod p$.
\item  If $K=\Q(\sqrt d)$ is a quadratic field, $\mathcal{P}$ is a prime ideal of
		$\mathcal{O}_K$ above a prime number $p$ and $\alpha\in \mathcal{O}_K$, then according to \cite{Lm00}, we have : $$\displaystyle\left(\dfrac{\alpha}{\mathcal{P}}\right)=\left\{
		\begin{array}{ll}
		\left(\dfrac{N(\alpha)}{p}\right) & \hbox{if $p$ is inert in $K$;} \\
		\left(\dfrac{\alpha}{p}\right)^{f_{\mathcal{P}}} & \hbox{if $\alpha\in\Q$;} \\
			\left(\dfrac{\alpha}{\mathcal{P}'}\right)	\left(\dfrac{N(\alpha)}{p}\right)  &  \hbox{if} \   p\mathcal{O}_K=\mathcal{P}\mathcal{P}', \mathcal{P}\neq \mathcal{P}'. \\
		\end{array}
		\right.
		$$
\end{enumerate}	
	\end{rem}
{\bf II. The Hilbert Symbol}.\\
Let $\mathcal{P}$ be  a prime  of $K$ (finite or infinite) and let  $K_\mathcal{P}$ be the completion of $K$ at $\mathcal{P}$.  For   $\alpha, \beta\in K_\mathcal{P}$, we define the local Hilbert Symbol as
$$
(\alpha, \beta)_\mathcal{P}=\begin{cases}1&\mbox{ if }H_{K_\mathcal{P}}(\alpha, \beta)\mbox{ splits };\\-1&\mbox{ if  not.}\end{cases}
$$
Recall that the Hilbert symbol
 $(a, b)_{\mathcal{P}}$ is equal to $-1$ if and only if  the equation $ax^2+by^2=1$ has no solutions in $K_\mathcal{P}$.\\
If $i_{\mathcal{P}}$ is the natural injection from $K$ to $K_\mathcal{P}$, then for   $\alpha, \beta\in K$ we define the global Hilbert Symbol as
$$
\left(\dfrac{\alpha, \beta}{\mathcal{P}}\right)=i_{\mathcal{P}}^{-1}((i_{\mathcal{P}}(\alpha), i_{\mathcal{P}}(\beta))_\mathcal{P}).
$$
\begin{thm}[\cite{Gr03}]
	The Hilbert symbol satisfies the following conditions:
\begin{enumerate}[\rm\indent i.]
	\item  $\left(\dfrac{\alpha,
		\beta}{\mathcal{P}}\right)=\left(\dfrac{\beta,
		\alpha}{\mathcal{P}}\right)$;
	\item \label{1:018}If $\mathcal{P}$ is a prime ideal  of $K$ unramified  in $K(\sqrt\alpha)$, then $\left(\dfrac{\alpha, \beta}{\mathcal{P}}\right)=\displaystyle\left(\dfrac{\alpha}{\mathcal{P}}\right)^{v_{\mathcal{P}}(\beta)}$;
	\item
	if   $\mathcal{P}_{\infty}$ denotes an infinite prime, then
	$\left(\dfrac{\alpha,
		\beta}{\mathcal{P}_{\infty}}\right)=-1$ if and only if $i_{\mathcal{P}}(\alpha)<0$ and $i_{\mathcal{P}}(\beta)<0$;
	\item $\beta$ is a norm in $K(\sqrt{\alpha})$ if and only if $\left(\dfrac{\alpha,
		\beta}{\mathcal{P}}\right)=1$ for all prime $\mathcal{P}$;
	\item
	$\displaystyle\prod_{\mathcal{P}}\left(\dfrac{\alpha,
		\beta}{\mathcal{P}}\right)=1$    $($the product formula$)$.
\end{enumerate}
\end{thm}
\begin{rem}\label{ghent2}
	Let $K$ be a number field and  $\alpha$, $\beta\in K^*=K\backslash\{0\}$. Then $H_K(\alpha, \beta)$ splits if and only if $\mathrm{Ram}(H_K(\alpha, \beta))=\emptyset$.  More generally, we have $$ \mathrm{Ram}(H_K(\alpha, \beta))=\left\{\mathcal{P} \text{ the prime ideal of } K  \;:\; \left(\dfrac{\alpha,
		\beta}{\mathcal{P}}\right)=-1\right\}.$$
Note also that the following splitting criterion for a quaternion algebras is well  known \cite[Corollary 1.10]{alsina}:
  the quaternion algebra $H_K(\alpha, \beta)$  is split if and only if
 its discriminant
  $D_{H_K(\alpha, \beta)}$ is equal to  the ring of integers $ \mathcal{O}_{K}$ of $K$.
\end{rem}
\section{Main Results}
Let $K=\mathbb{Q}(\sqrt{d})$ be a quadratic number field and $\mathcal{O}_K$ its ring of integers, where $d\neq1$ is a squarefree integer.
\subsection{First case: $\alpha\in \mathcal{O}_K$   and  $m\equiv1\pmod4$ a squarefree integer}
 Let us begin  by defining an algebraic integer $\alpha\in \mathcal{O}_K$  to be odd if any  prime ideal appearing in its factorization into prime ideals in $\mathcal{O}_K$ does not lie above $2$.

 The numbers  $m$ and $\alpha$ are said to satisfy   hypotheses (H) if they satisfy the following two conditions:
\begin{enumerate}[\indent\rm1.]
	\item $\alpha$ is an odd integer of $\mathcal{O}_K$, \ $m\equiv 1\pmod 4$ and $d$  are two distinct  squarefree  integers with $m>0$;
	\item $\alpha$ and $m$ are relatively prime.
\end{enumerate}

\noindent Moreover, we adopt, in the sequel, the following  notation:
\begin{enumerate}[\indent\rm1.]
	\item $\dfrac{m}{\gcd(m,d)}=(q_1q_2\cdots q_t)^\lambda(r_1r_2\cdots r_s)^\gamma$ where $\lambda,\ \gamma\in\{0, 1\} $ and $q_i$, $r_i$ are prime numbers satisfying  $\left(\dfrac{d}{q_i}\right)=-1$ if $\lambda=1$, and $\left(\dfrac{d}{r_i}\right)=1$ if $\gamma=1$ for all $i$;
	\item $\mathcal{Q}_i$ (resp. $\mathcal{R}_i$, $\mathcal{R}^{'}_i$) denotes the prime ideal (resp. ideals) of $\mathcal{O}_K$ above $q_i$ (resp. $r_i$);
	\item   $\mathcal{I}$ denotes  the set of all prime ideals $\mathcal{P}$ of $\mathcal{O}_K$ such that $v_{\mathcal{P}}(\alpha)$ is odd and
	whose absolute norm is a rational prime.
\end{enumerate}
We can now establish our first main result.
\begin{thm}\label{theorem4} If $m$ and $\alpha$ satisfy  the hypotheses $(\mathrm{H})$, then $H_K(\alpha, m)$ splits if and only if  the following conditions are satisfied:
\begin{enumerate}[\indent\rm1.]
	\item $\left(\dfrac{N(\alpha)}{q_i}\right)= \left(\dfrac{\alpha}{\mathcal{Q}_i}\right)=1$, if $\lambda=1$;
	\item $\left(\dfrac{N(\alpha)}{r_j}\right)=\left(\dfrac{\alpha}{\mathcal{R}_j}\right)=1$, if $\gamma=1$;
	\item for any prime $\mathcal{P}\in\mathcal{I}$ we have $\left(\dfrac{m}{N(\mathcal{P})}\right)=1$.
\end{enumerate}
\end{thm}
\begin{proof}
 Let $\mathcal{P}$ be  an odd prime ideal of $K$. Lemma 1 of \cite{B-P-76} implies that the discriminant $D(L/K)=(\mathcal{Q}_1\mathcal{Q}_2\cdots \mathcal{Q}_t)^\lambda(\mathcal{R}_1\mathcal{R}^{'}_1\mathcal{R}_2\mathcal{R}^{'}_2\cdots \mathcal{R}_s\mathcal{R}^{'}_s)^\gamma$, where $L=K(\sqrt m)$.\\
- If $\mathcal{P}\neq \mathcal{Q}_i$, $\mathcal{P}\neq \mathcal{R}_j$ and $\mathcal{P}\neq \mathcal{R^{'}}_j$, then $\mathcal{P}$ is unramified in $K(\sqrt{m})$. Hence by putting  $N(\mathcal{P})=p^{f_{\mathcal{P}}}$, Remark \ref{ghent1} implies that
$$\left(\dfrac{\alpha, m}{\mathcal{P}}\right)=\displaystyle\left(\dfrac{m}{\mathcal{P}}\right)^{v_{\mathcal{P}}(\alpha)}=\left(\dfrac{m}{p}\right)^{f_{\mathcal{P}}v_{\mathcal{P}}(\alpha)}=\begin{cases}\left(\dfrac{m}{p}\right),&\mbox{ if }\mathcal{P}\in\mathcal{I};\\1,&\mbox{ if  not.}\end{cases}$$
- If $\mathcal{P}=\mathcal{Q}_i$, $\mathcal{R}_j$ or $\mathcal{R^{'}}_j$, then, since $\alpha$ and $m$ are relatively prime,  $v_{\mathcal{P}}(\alpha)=0$ and $v_{\mathcal{P}}(m)=1$. Hence $\mathcal{P}$ is unramified in $K(\sqrt{\alpha})$ and thus
$$\left(\dfrac{\alpha, m}{\mathcal{P}}\right)=\displaystyle\left(\dfrac{\alpha}{\mathcal{P}}\right)^{v_{\mathcal{P}}(m)}=\begin{cases}\left(\dfrac{N(\alpha)}{q_i}\right),&\mbox{ if }\mathcal{P}=\mathcal{Q}_i;\\\left(\dfrac{\alpha}{\mathcal{P}}\right),&\mbox{ if  not.}\end{cases}$$
 $\bullet$ As $m\equiv 1\pmod 4$, then the prime ideal $\mathcal{P}_2$ of $\mathcal{O}_K$ above $2$ is unramified in $K(\sqrt{m})$, and since $\alpha$ is an odd integer  (hence $v_{\mathcal{P}_2}(\alpha)=0$), we have:
$$\left(\dfrac{\alpha, m}{\mathcal{P}_2}\right)=\displaystyle\left(\dfrac{m}{\mathcal{P}_2}\right)^0=1.$$
$\bullet$ Now let $\mathcal{P}_\infty$ be an infinite prime ideal, then $\left(\dfrac{\alpha, m}{\mathcal{P}_\infty}\right)=1$, because $i_{\mathcal{P}_\infty}(m)=m>0$. Finally, Remarks \ref{ghent1} and \ref{ghent2} imply  the assertions.
\end{proof}
\begin{rem}
Thanks to   Hilbert's symbol product formula, we get the same results if we consider the conditions $m\equiv1\pmod4$ and $\alpha\in\mathcal{O}_K$ with $d\not\equiv 1\pmod 8$ instead of the conditions $m\equiv1\pmod4$  and $\alpha$ is an odd integer of $\mathcal{O}_K$.
In this case we will say that $m$ and $\alpha$ satisfy  the hypothesis $\mathrm{\hat{H}}$.
\end{rem}
The following corollary generalizes Theorem 3.6 and Proposition 3.7 of \cite{Sav16}.
\begin{cor}\label{2}
 If $m$ and $\alpha$ satisfy  hypotheses $\mathrm{\hat{H}}$, $\alpha\in\Z$ and $m$ is prime, then $H_K(\alpha, m)$ splits if and only if one of the following conditions holds.
\begin{enumerate}[\indent\rm1.]
	\item $\alpha$ is a square in $\N$.
	\item $-\alpha$ is a square in $\N$ and $m$ divides $d$ or $d$ is a quadratic residue modulo $m$.
	\item $\alpha =\pm \delta t^2$ with $\delta>1$ and each prime divisor of $\delta$, which splits in $K$, is a quadratic residue modulo $m$ and
	\begin{enumerate}[-]
		\item either $m$ divides $d$ or $d$ is  not a quadratic residue modulo $m$
		\item or $d$ and $\alpha$ are quadratic residue modulo $m$.
	\end{enumerate}
\end{enumerate}	

\end{cor}
\begin{cor}
If $\alpha$ is a unit of  $K$ and $m$ divides $d$, then  $H_K(\alpha, m)$ splits.
\end{cor}
\begin{proof}
It is enough to note  that if $\alpha$ is a unit of $K$, then $\mathcal{I}=\emptyset$; and  since $m$ divides  $d$ we have $\lambda=\gamma=0$.
\end{proof}
\subsection{Applications of  Theorem \ref{theorem4}}
In this subsection, we state some applications of the first main theorem. We start with a classical lemma \cite{KaWa}:
\begin{lem}Let $K$ be a quadratic field.
\begin{enumerate}[\rm1.]
	\item $K(\sqrt{\alpha})$ is a biquadratic number field if and only if $N(\alpha)$ is a square in $\N$.
	\item $K(\sqrt{\alpha})$ is a quartic cyclic number field if and only if $dN(\alpha)$ is a square in $\N$.
\end{enumerate}
\end{lem}
\noindent The first application of the main result is the following theorem.
\begin{thm}\label{ghent7}
If $m$ and $\alpha$ satisfy the  hypotheses $(\mathrm{H})$ and $K(\sqrt{\alpha})=\Q(\sqrt{d}, \sqrt{d_1})$ is a biquadratic number field, then $H_K(\alpha, m)$ splits if and only if  the following conditions are satisfied:
\begin{enumerate}[\rm1.]
	\item $\left(\dfrac{d_1}{r_j}\right)=1$, if $\gamma=1$.
	\item For any prime $p|d_1$ such that  $(p, d)=1$ and $\left(\dfrac{d}{p}\right)=1$, we have  $\left(\dfrac{m}{p}\right)=1$.
\end{enumerate}
\end{thm}
\begin{proof}
As $K(\sqrt{\alpha})=\Q(\sqrt{d}, \sqrt{d_1})$ is a biquadratic number field, then the previous Lemma imply the $N(\alpha)$ is a square; moreover  $\alpha =t^2d_1$ with $t\in \Q(\sqrt{d})$.  So  $H_K(\alpha, m)$ splits if and only if  the following conditions are satisfied:
\begin{enumerate}[\indent\rm1.]
	\item $\left(\dfrac{\alpha}{\mathcal{R}_j}\right)=\left(\dfrac{t^2d_1}{\mathcal{R}_j}\right)=\left(\dfrac{d_1}{r_j}\right)=1$, if $\gamma=1$.
	\item For any prime $\mathcal{P}\in\mathcal{I}$, we have $\left(\dfrac{m}{N(\mathcal{P})}\right)=1$.
\end{enumerate}
If a prime $\mathcal{P}\in\mathcal{I}$, then $v_{\mathcal{P}}(\alpha)=2v_{\mathcal{P}}(t)+v_{\mathcal{P}}(d_1)=2v_{\mathcal{P}}(t)+v_{p}(d_1)v_{\mathcal{P}}(p)$ is odd and $N(\mathcal{P})=p$ is a prime integer, which is equivalent to  $p|d_1$, $(p, d)=1$,  $\left(\dfrac{d}{p}\right)=1$ and $\left(\dfrac{m}{p}\right)=1$.
\end{proof}
\noindent As a second   application of the main Theorem, we assume that  $K(\sqrt{\alpha})$ is  a quartic cyclic number field.
\begin{thm}
	If $m$ and $\alpha$ satisfy  hypotheses $(\mathrm{H})$, $d\equiv 1\pmod 4$ and $K(\sqrt{\alpha})=\Q(\sqrt{ a+b\sqrt{d}})$ is quartic cyclic number field $(a^2=d(b^2+c^2))$, then $H_K(\alpha, m)$ splits if and only if  the following conditions are satisfied:
	\begin{enumerate}[\rm1.]
		\item  $\lambda=0$ and $(-1)^{(r_j-1)(d-1)/8}\left(\dfrac{2}{r_j}\right)^b\left(\dfrac{r_j}{d}\right)_4=1$,  if $\gamma=1$,
		\item For any prime $\mathcal{P}\in\mathcal{I}$, we have $\left(\dfrac{m}{N(\mathcal{P})}\right)=1$.
	\end{enumerate}
\end{thm}
\begin{proof}
It is an immediate consequence of  the main theorem of \cite{WiHaFr85}.
\end{proof}
\noindent As a third   application of the main Theorem,
 we give necessary and sufficient conditions for $H_K(\alpha, m)$ to split over $K$ whenever  $\alpha=\varepsilon_d$  is a fundamental unit of $K=\Q(\sqrt{d})$.
\begin{thm}\label{ghent4}
	Let's assume that the norm of $\varepsilon_d$ is $-1$, and let $d = p_1\cdots p_n$. Suppose that  $p_j$ and $m$ are primes
	$\equiv 1\pmod 4$ such that $\gamma=1$ and $\left(\dfrac{p_j}{r_i}\right) = 1$ for all $1 \leq j \leq n$ and $1 \leq i \leq s$. Then $H_K(\varepsilon_d, m)$ splits if and only if $$r_i\equiv 1\pmod 4\text{ and }\left(\dfrac{r_i}{d}\right)_4\left(\dfrac{d}{r_i}\right)_4=1,$$
for all $1 \leq i \leq s$.	
\end{thm}
\begin{proof} Note that $m$ and $\alpha=\varepsilon_d$  satisfy the hypotheses (H), $\mathcal{I}=\emptyset$. Then  \cite[Corollary on page 143]{Fur}
gives the  claimed result.
\end{proof}
 We now suppose that the norm of $\varepsilon_d$ is $1$. It is known that $\varepsilon_d$ does not always belong to  $\Z[\sqrt{d}]$. Actually,
 $\varepsilon_d\in\Z[\sqrt{d}]$ if $d\equiv 3, 2\pmod 4$ or $d\equiv 1\pmod 8$, but if $d\equiv 5\pmod 8$, $\varepsilon_d$ may not belong to $\Z[\sqrt{d}]$. To simplify,
  we set $\varepsilon=\varepsilon_d^3$ if $d\equiv 5\pmod 8$ and $\varepsilon_d\notin\Z[\sqrt{d}]$, and
  $\varepsilon=\varepsilon_d$ otherwise.
\begin{thm}\label{ghent6}
Let	$m$ and $d$ be   squarefree  integers such that $\gamma=1$, let the norm of $\varepsilon=x+y\sqrt{d}$ be  $1$.
 Then $H_K(\varepsilon_d, m)$ splits if and only if one of the following conditions holds.
\begin{enumerate}[\indent\rm1.]
	\item $\left(\dfrac{d_1}{r_j}\right)=1$ where  $d_1$ divides  $d$ and $2d_1(x+1)$ is a square in $\N$ for all $1 \leq i \leq s$.
	\item $\left(\dfrac{2}{r_j}\right)=\left(\dfrac{d_1}{r_j}\right)$ where  $d_1$ divides  $d$ and $d_1(x+ 1)$ is a square in $\N$ for all $1 \leq i \leq s$.
\end{enumerate}
\end{thm}
\begin{proof}
As the norm of $\varepsilon_d$ is $1$, then $\varepsilon$ is too, i.e., $(x+1)(x-1)=dy^2$. Furthermore $K(\sqrt{\varepsilon})$ is a biquadratic field.\\
$\bullet$ If $x$ is odd, then according to the unique prime factorization
 in $\Z$, there exist two  squarefree  integers $d_1$ and $d_2$ such that  $d=d_1d_2$, $(x+1)=2d_1y_1^2$ and  $(x-1)=2d_2y_2^2$
 with $y_1, y_2\in\N$ and $2y_1y_2=y$. In this case, one can easily check that $2d_1(x+1)$ is a square in
  $\N$ and that $d_1\varepsilon=(y_1d_1+y_2\sqrt{d})^2$. With Theorem \ref{ghent7}, we have
$H_K(\varepsilon_d, m)$ splits if and only if $\left(\dfrac{d_1}{r_j}\right)=1$.\\
$\bullet$ If $x$ is even, then proceeding similarly, we get that $d_1(x+ 1)$
is a square in $\N$ and   $2d_1\varepsilon=(y_1d_1+y_2\sqrt{d})^2$. 
As in the first case  $H_K(\varepsilon_d, m)$ splits if and only if $\left(\dfrac{2}{r_j}\right)=\left(\dfrac{d_1}{r_j}\right)$.
\end{proof}
\begin{rem}
	Note that if the norm of  $\varepsilon_d$ is $1$, then $d_1$ exists and it is unique.
\end{rem}
\subsection{Second case: $m$ and $\alpha$ are prime integers}
In this subsection, we replace $\{\alpha, m\}$ by $\{p, q\}$, where $p$ and $q$ are rational primes. The purpose  is then to establish the following theorem that classify division quaternion algebras  $H_{\Q(\sqrt d)}(p, q)$ over quadratic fields where $d\neq1$ is a squarefree integer:
\begin{thm}
\label{maintheorem}
{Let} $d\neq1$  {be a squarefree integer} and let
 $K=\mathbb{Q}(\sqrt{d})$,
 with discriminant  $\Delta_{K}$. {Let} $p$ and $q$  be two  positive primes.
 {Then the quaternion algebra} $H_{K}(p,q)$   {is a division algebra if and only if} one of the following conditions
  holds:
 \begin{enumerate}[\rm1.]
 \item
$p$ and $q$  are odd and distinct, $p$ or $q \equiv 1 \pmod4$,
$(\frac{p}{q})=-1$
     and
    $(\frac{\Delta_{K}}{p})=1$   {or}  $(\frac{\Delta_{K}}{q})=1$;
 \item
  $q=2$,  $p \equiv 3\pmod8$   or      $p \equiv 5\pmod8$ and
 either $(\frac{\Delta_{K}}{p})=1$   {or}
  $d\equiv1\pmod8$;
 \item
 $p$ and $q$ are odd,  with $p\equiv q\equiv 3\pmod4$, and
 \begin{itemize}
 \item
  $(\frac{q}{p})\neq 1$ and
 either $(\frac{\Delta_{K}}{p})=1$   {or} $d$$\equiv$$1$ $\pmod8$;\\
 or
 \item
  $(\frac{p}{q})\neq 1$ and
 either $(\frac{\Delta_{K}}{q})=1$   {or} $d$$\equiv$$1$ $\pmod8$.
 \end{itemize}
 \end{enumerate}
\end{thm}
\begin{proof}
The proof of the first assertion is a simple deduction from Theorem \ref{ghent7}. It suffices to take $\alpha, m\in\{p, q\}$ such that $\alpha\neq m$ and $m\equiv1\pmod4$. Assertions 2 and 3 follow easily applying other results which we obtained in this article, namely Propositions \ref{threedotseven}, \ref{threedoteight} and \ref{threedotfourteen} below, after taking into account
Lemma \ref{threedotfour} and Lemma \ref{threedotthirteen} below.
\end{proof}
To complete the proof of this second main result, we need some preliminary results. Let $K$ be a number field  and $\mathcal{O}_{K}$  its ring of integers.
If $\mathcal{O}_{K}$ is a principal ideal domain, then we may identify the ideals
of $\mathcal{O}_{K}$ with their generators, up to units.
Thus, in a quaternion algebra $H$ over $\mathbb{Q}$, the element $D_{H}$ (the discriminant of $H$) turns out to be an integer, and $H$  is split if and only if $D_{H}=1$.
On the other hand,  a quaternion algebra $H_\Q (\alpha, \beta)$  is a division algebra if and only if there is   a prime
$p$ such that $p | D_{H_\Q  (\alpha, \beta)}$.
We continue  this note with two statements following from the classical
Albert-Brauer-Hasse-Noether theorem. Proofs of specific formulations
of this theorem can be found in \cite{linowitz, chinburg}.
\begin{thm}
\label{twodotnine}
Let $H_{F}$   be a quaternion algebra over a number field $F$   {and let} $K$   {be a quadratic   extension of} $F.$   {Then there is an embedding of} $K$   {into}
$H_{F}$   {if and only if no prime of} $F$   {which ramifies in} $H_{F}$   splits in  $K.$
\end{thm}
\begin{prop} \label{twodotten}
 Let $F$   be a number field and let  $K$   be a quadratic extension of $F$.   Let   $H_{F}$   be a quaternion algebra over  $F$.   Then   $K$   splits $H_{F}$   if and only if there exists an $F$-embedding   $K\hookrightarrow H_{F}$.
\end{prop}
In \cite{Sav16} the second author obtained the following result about quaternion algebras over the field $\mathbb{Q}(i)$ which is also a simple deduction from Corollary \ref{2}.
\begin{prop}
\label{threedotone}
Let $p\equiv 1$ $\pmod4$    be a prime
integer and  let $m$  be an integer which is not a quadratic
residue modulo $p$. Then the quaternion algebra $H_{\mathbb{Q}( i) }( m ,p) $  is a division
algebra.
\end{prop}
In    \cite{Sav17} the second author obtained some sufficient conditions for a quaternion algebra
$H_{\mathbb{Q}( i) }        (p, q)$ to split, where $p$ and $q$ are two distinct   primes:
\begin{prop}
\label{threedottwo}
Let $d\neq 0, 1$   {be a   squarefree integer} such that $d\not\equiv 1$ $\pmod8$,   {and let} $p$ and $q$   {be two primes,} with $q\geq 3$ and $p\neq q.$   {Let} $\mathcal{O}_{K}$   {be the ring of integers of the quadratic field} $K=\mathbb{Q}(\sqrt{d})$   {and} let $\Delta_{K}$   be its discriminant.
\begin{enumerate}[\rm1.]
\item
    {If} $p\geq 3$  and  both   $(\frac{\Delta_{K}}{p} ) $ and  $(\frac{\Delta_{K}}{q})$ are not equal to $1$,
     {then the quaternion algebra} $H_{\mathbb{Q}(\sqrt{d})}(p,q)$   {splits};
\item   {If} $p=2$   and   $(\frac{\Delta_{K}}{q})\neq 1,$   {then the quaternion algebra} $H_{\mathbb{Q}(\sqrt{d})}(2, q)$   {splits}.
\end{enumerate}
\end{prop}
From the aforementioned results we   deduce easily  a necessary and sufficient condition for a quaternion algebra $H_{\mathbb{Q}(i)}(p, q)$ to be a division algebra:
\begin{prop}
\label{threedotthree}
  {Let} $p$ and $ q$   {be two distinct odd primes,}  {such that  } $(\frac{q}{p})\neq 1$.  {Then the quaternion algebra} $H_{\mathbb{Q}(i)}(p,q)$   {is a division algebra if and only if} $p\equiv 1\pmod 4$   {or} $q\equiv 1\pmod 4$.
\end{prop}
\begin{proof}
To prove the necessity,  note that if $H_{\mathbb{Q}(i)}(p,q)$ is a division algebra, then   Proposition \ref{threedottwo} and Theorem \ref{twodoteightone} tell us that
  $(\frac{\Delta_{K}}{p})=1$ or $(\frac{\Delta_{K}}{q})= 1$. This is equivalent to $p\equiv 1$ (mod $4$)   {or} $q\equiv 1$ (mod $4$).

To prove the sufficiency, we must distinguish amongst two cases:
\begin{itemize}
\item
$p\equiv 1$ (mod $4$): \\
Since $(\frac{q}{p})\neq 1$,    Proposition \ref{threedotone}  tells us that $H_{\mathbb{Q}(i)}(p,q)$ is a division algebra.
\item
$q\equiv 1$ (mod $4$): \\
Since $(\frac{q}{p})=-1$, the quadratic reciprocity law implies  $(\frac{p}{q})=-1$.
   Proposition \ref{threedotone}  then tells us that   $H_{\mathbb{Q}(i)}(p,q)$ is a division algebra.
\end{itemize}
\end{proof}
We ask ourselves whether we can obtain a necessary and sufficient explicit condition for
$H_{\mathbb{Q}(\sqrt{d})}(p,q)$ to be a division algebra when $d$ is arbitrary.
From Proposition \ref{threedottwo} we obtain a necessary explicit condition for   $H_{\mathbb{Q}(\sqrt{d})}(p,q)$ to be a division algebra, namely:   {if} $H_{\mathbb{Q}(\sqrt{d})}(p,q)$   {is a division algebra, then} $(\frac{\Delta_{K}}{p})=1$   {or} $(\frac{\Delta_{K}}{q})=1$. However  this condition is not sufficient:
for example, if we let $K=\mathbb{Q}(\sqrt{3}),\ p=7, \ q=47$,
then $(\frac{\Delta_{K}}{p})\neq 1$ and $(\frac{\Delta_{K}}{q})=1$,
the quaternion algebra $H_{\mathbb{Q}}(7,47)$ is a division algebra, but
the quaternion algebra $H_{\mathbb{Q}(\sqrt{3})}(7,47)$ is not a division algebra.

It is known   \cite{kohel} that  if a prime integer $p$ divides  $D_{H(\alpha, \beta)}$, then it must divide  $2\alpha\beta$,
hence we may restrict our attention to these primes.
In other words, in order to obtain a sufficient condition for a quaternion algebra
$H_{\mathbb{Q}(\sqrt{d})}(p,q)$   to be a division algebra,
 it is important to study the ramification of the primes   $ 2,p,q $     in the   algebra    $H_\Q(p,q)$.
The following lemma  \cite[Lemma 1.21]{alsina} gives us a hint:
\begin{lem}
\label{threedotfour}
Let  $p$ and $q$ be two primes, and let $H_{\mathbb{Q}}(p,q)$ be a quaternion algebra of discriminant $D_H$.
\begin{enumerate}[\rm1.]
\item
if $p\equiv q\equiv 3$ $\pmod4$   {and} $(\frac{q}{p})\neq 1$,   {then} $D_{H}=2p$;
 \item
if $q=2$ and  $p\equiv 3$ $\pmod8$,  {then} $D_{H}=pq=2p$;
\item
if  $p$   {or} $q\equiv1 \pmod 4$, with  $p\neq q$  {and} $(\frac{p}{q})=-1$,   then  $D_{H}=pq$.
\end{enumerate}
\end{lem}
In addition, the following lemma  \cite[Lemma 1.20]{alsina}     tells us
precisely  when a quaternion algebra $H_{\mathbb{Q}}(p,q)$ splits.
\begin{lem}
\label{threedotfive}
Let  $p$ and $q$ be two prime integers. Then
$H_{\mathbb{Q}}(p,q)$ is a matrix algebra if and only if one of
the following conditions is satisfied:
\begin{enumerate}[\rm1.]
\item
$p=q=2$;
 \item
$p=q \equiv1 \pmod 4$;
\item
$q=2$   {and} $p\equiv\pm 1$ $\pmod8$;
\item
$p\neq q,$ $p\neq 2,$ $q\neq 2,$ $(\frac{q}{p})=1$ and either $p$ or $q$ is congruent to $1\pmod4$.
\end{enumerate}
\end{lem}
The next theorem  \cite[Theorem 1.22]{alsina}  describes  the discriminant of
$H_{\mathbb{Q}}(p,q)$, where $p$ and $q$ are primes:
\begin{thm}
\label{teo}
Let $H = ( \frac{a,b}{\Q} )$ be a quaternion algebra.
\begin{enumerate}[\rm1.]
\item
If $D_H=1$, then $H$ splits;
\item
If $D_H=2p$, $p$ prime and $p \equiv 3 \pmod{ 4}$, then $H \cong ( \frac{p, -1}{\Q} )$;
\item
If $D_H=pq$, $p,q$ primes, $q \equiv 1 \pmod{4}$ and $(\frac{p}{q})=-1$, then $H \cong ( \frac{p,\; q}{\Q} )$.
\end{enumerate}
If $a$ and $b$ are prime numbers, then the algebra $H$ satisfies one and only one of the above statements.
\end{thm}

 Note that when  $q=2$ and $p$ is a prime such that $p \equiv 3$ (mod 8),  then,
  according to Lemma \ref{threedotfour} the discriminant
$D_{H_{\mathbb{Q}(p,q)} }$ is equal to $2p$,
so $H_{\mathbb{Q}}  (p,q  )$ is a division algebra. The
next proposition shows what   happens when we extend the field of scalars from $\Q$ to $\Q(\sqrt{d})$.
\begin{prop}
\label{threedotseven}
 {Let} $p$   {be an odd prime,} with $p \equiv 3$ $\pmod8$.   Let  $K=\mathbb{Q}(\sqrt{d})$   {and let} $\Delta_{K}$   {be the discriminant of} $K.$
  Then  $H_{\mathbb{Q}(\sqrt{d})}(p,2)$   {is a division algebra if and only if} $(\frac{\Delta_{K}}{p})= 1$   {or} $d$ $\equiv 1$ $\pmod8$.
\end{prop}
\begin{proof}
 If $H_{\mathbb{Q}(\sqrt{d})}(p,2)$ is a division algebra then, from
  Proposition \ref{threedottwo}, Theorem \ref{twodoteightone},
 Theorem \ref{twodotnine} and Proposition \ref{twodotten}, we conclude that $(\frac{\Delta_{K}}{p})=1$ or $d$ $\equiv 1$ (mod $8$).

Conversely, since $p \equiv 3$ (mod $8$) then, according to Lemma \ref{threedotfour}(ii) we must have $D_{H}=2p.$
 It follows that the primes which ramify in   $H_{\mathbb{Q}}(p,2)$ are precisely
$p$ and $2$. Since either $(\frac{\Delta_{K}}{p})= 1$ or $d \equiv 1$ (mod $8$) then,
 by the decomposition of primes in quadratic fields, we obtain that either $p$ or $2$ splits in the ring or integers of $K$.
From Theorem \ref{twodotnine}, Proposition  \ref{twodotten} and Theorem \ref{twodoteightone},
we conclude that  $H_{\mathbb{Q}(\sqrt{d})}(p,2)$ is a division algebra.
\end{proof}
We next study the case where $p$ and $q$ are primes, both congruent to $3$ modulo $4$.
If $(\frac{q}{p})\neq 1$, then,
according to Lemma \ref{threedotfour}(i), the discriminant
$D_{H_{\mathbb{Q}}   (p,q)}$ is equal to $2p$, so $H_{\mathbb{Q}}(p,q)$ is a division algebra.
The next proposition tells us when
the quaternion algebra $H_{\mathbb{Q}(\sqrt{d})}(p,q)$ is still a division algebra.
\begin{prop}
\label{threedoteight}
 {Let} $p$ and $q$   {be two odd prime integers,} with $p\equiv q\equiv 3$ ({mod} $4$)  and  $(\frac{q}{p})\neq 1$.
 Let  $K=\mathbb{Q}(\sqrt{d})$   {and let} $\Delta_{K}$   {be the discriminant of} $K.$
  {Then the quaternion algebra} $H_{\mathbb{Q}(\sqrt{d})}(p,q)$   {is a division algebra if and only if} $(\frac{\Delta_{K}}{p})= 1$   {or} $d$ $\equiv 1$ $\pmod8$.
\end{prop}
\begin{proof}
  If $H_{\mathbb{Q}(\sqrt{d})}(p,q)$ is a division algebra then from Proposition \ref{threedottwo}(i) it follows that either
   $(\frac{\Delta_{K}}{p})=1$ or $(\frac{\Delta_{K}}{q})=1$. But, according to Lemma \ref{threedotfour}(i) we must have
    $D_{H_{\mathbb{Q}}(p,q)}=2p.$ So the integral primes   which
ramify in $H_{\mathbb{Q}}(p,q)$   and could split in $K$ are precisely $p$ and $2.$
Finally, after applying Theorem \ref{twodoteightone}, Theorem \ref{twodotnine}, Proposition \ref{twodotten} and  the decomposition of primes in quadratic fields,
we obtain that either $(\frac{\Delta_{K}}{p})=1$  or $d$ $\equiv 1$ (mod $8$).
The proof of the converse is similar to  the proof  of sufficiency of Proposition \ref{threedotseven}.
\end{proof}

The only case left  out is $q = 2, \  p\equiv 5$ $\pmod8$.
We consider first the quaternion algebra
$H_{\mathbb{Q}}(p,q)$, and we get  the following result:
\begin{lem}
\label{threedotthirteen}
Let  $p\equiv 5$ $\pmod8$ be a prime integer.
 Then the discriminant of the quaternion algebra $H_{\mathbb{Q}}(p,2 )$
is equal to $2p$, and hence $H_{\mathbb{Q}}(p,2)$ is a division algebra.
\end{lem}
\begin{proof}
We give here a simple proof which is independent of the theorems stated above.
We know that if a prime divides the discriminant of $H_{\Q}(a,b)$
then it must divide $2ab$.
Since $p\equiv 5$ (mod $8$), from the properties of the Hilbert symbol   and of the
Legendre symbol we obtain:
$
(2,p)_{p}=(\frac{2}{p})=(-1)^{\frac{p^{2}-1}{8}}=-1
$
and
$
(2,p)_{2}=(-1)^{\frac{p-1}{2}\cdot \frac{1-1}{2} + \frac{p^{2}-1}{8}}=-1.
$
Hence the primes which ramify in $H_{\mathbb{Q}}(p,2)$ are exactly
$2$ and $p$. Therefore,   the reduced discriminant of  $H_{\mathbb{Q}}(p,2 )$  must be equal to $2p$.
\end{proof}
We turn now our attention   to the quaternion algebra $H_{\mathbb{Q}(\sqrt{d})}(p,q)$, where $q = 2$ and $p\equiv 5$ $\pmod8$.
\begin{prop}
\label{threedotfourteen}
 {Let} $p$   {be an odd prime,} with $p \equiv 5$ $\pmod8$.   Let  $K=\mathbb{Q}(\sqrt{d})$   {and let} $\Delta_{K}$   {be the discriminant of} $K.$ Then  $H_{\mathbb{Q}(\sqrt{d})}(p,2)$   {is a division algebra if and only if} $(\frac{\Delta_{K}}{p})= 1$   {or} $d$ $\equiv 1$ $\pmod8$.
\end{prop}
\begin{proof}
The proof is similar to the proof of Proposition \ref{threedotseven},
after replacing   Lemma \ref{threedotfour} with   Lemma \ref{threedotthirteen}.
\end{proof}

Taking into account these results and  Theorems \ref{twodoteightone} and \ref{maintheorem}, we are able to understand when a quaternion algebra $H_{\mathbb{Q}(\sqrt{d})}(p,q)$ splits. It is clear also that in each one of the cases  covered by Lemma \ref{threedotfive},
 a quaternion algebra $H_{\mathbb{Q}(\sqrt{d})}(p,q)$ splits.

\section{Examples involving  prime Fibonacci numbers}
Let $(F_n)_{n\geq0}$ be the   Fibonacci sequence which  is  defined by the following recurrence relation:
\begin{center} $F_0=0$, $F_1=1$ and $F_n = F_{n-1} + F_{n-2}$ for $n\geq 2$.\end{center}
 {\it Binet's formula},  which discovered by the  mathematician J. P. Marie Binet (1786-
1856), states that:
\begin{equation*}
F_n=\dfrac{1}{\sqrt{5}}\left[\left(\dfrac{1 + \sqrt{5}}{2}\right)^{n}-\left(\dfrac{1 - \sqrt{5}}{2}\right)^{n}\right].
\end{equation*}
By means of this formula, one can   show that (see \cite{GrRa12}):
\begin{equation*}
F_{2n+1}=F_{n+1}^2 + F_{n}^2  \text{ and } 4^nF_{2n+1}=\sum_{k=0}^{n}5^k\mathrm{C}_{2n+1}^{2k+1}.
\end{equation*}
where $C_m$ denotes the $m$-th Catalan number.
\noindent From this formula
we deduce  that if $p$ is an odd prime number $>5$, then $F_{p}\equiv 1\pmod 4$ and $F_{p} \equiv \left(\dfrac{p}{5}\right) \pmod p$. This implies the following lemma.
\begin{lem}\label{ghent3}
	Let $p$ be a prime $>5$. Then the Fibonacci number $F_p$ has the following properties:
	\begin{enumerate}[\rm1.]
		\item If $p\equiv 1 \pmod 4$, then $\left(\dfrac{F_p}{p}\right)=1$,
		\item If $p\equiv 3 \pmod 4$, then $\left(\dfrac{F_p}{p}\right)=\left(\dfrac{p}{5}\right)$.
	\end{enumerate}
\end{lem}
To prove the last theorem of this paper, we also need the following lemma.
\begin{lem}[\cite{Ta}]\label{ghent5}
	Let $F_p$ be a Fibonacci number with prime index $p\equiv 1\pmod 4$. Then
	$$
	\left(\dfrac{F_p}{p}\right)_4\left(\dfrac{p}{F_p}\right)_4=\left\{
	\begin{array}{ll}
	1 & \hbox{if $p\equiv 1\pmod 3$,} \\
	\left(\dfrac{2}{p}\right) & \hbox{if not.}
	\end{array}
	\right.
	$$
\end{lem}
\begin{thm}
Let $F_p$ be a Fibonacci prime number  with $p>5$ and let $K=\Q(\sqrt{p})$ be a quadratic  number field. Then $H_K(\varepsilon_p, F_p)$ splits if and only if one of the following conditions holds.
\begin{enumerate}[\indent\rm1.]
	\item $p\equiv 1\pmod{12}$;
	\item $p\equiv 17\pmod{24}$;
	\item $p\equiv 3\pmod 4$, $p^2\equiv 1\pmod 5$ and  $\left(\dfrac{2}{F_p}\right)=1$;
	\item $p\equiv 3\pmod 4$ and  $p^2\equiv 4\pmod 5$.
\end{enumerate}
\end{thm}
\begin{proof}
First of all, recall that if $F_n$ is prime, then $n$ is prime;  and   $n$ divides $m$ if and only if $F_n$ divides $F_m$, hence since $F_p$ is a prime Fibonacci number such that $p> 5$, it follows that $p$ is an odd prime.\\
$\bullet$ If $p\equiv 1\pmod 4$, then $N(\varepsilon_p)=-1$ and $\left(\dfrac{p}{F_p}\right)=\left(\dfrac{F_p}{p}\right)=1$ ((1) of Lemma \ref{ghent3}). We can then apply Corollary \ref{ghent4} and find that $H_K(\varepsilon_p, F_p)$ splits if and only if $\left(\dfrac{F_p}{p}\right)_4\left(\dfrac{p}{F_p}\right)_4=1$.  According to Lemma \ref{ghent5}, this is equivalent to $p\equiv 1\pmod{3}$ or $p\equiv 1\pmod{8}$. By the Chinese remainder theorem, we deduce that $p\equiv 1\pmod{12}$ or $p\equiv 17\pmod{24}$.\\
$\bullet$ If $p\equiv 3\pmod 4$, then $N(\varepsilon_p)=1$. We must now study the value of the symbol $\left(\dfrac{F_p}{p}\right)$. If $\left(\dfrac{F_p}{p}\right)=1$, then the second assertion of   Lemma \ref{ghent3} implies that $\left(\dfrac{5}{p}\right)=1$, i.e., $p^2\equiv 1\pmod 5$. In this case $H_K(\varepsilon_p, F_p)$ splits if and only if $\left(\dfrac{2}{F_p}\right)=1$ (\cite[Lemma 6]{Az-00} and 2. of  Theorem \ref{ghent6}). Finally,  if $p^2\equiv 4\pmod 5$, then Theorem \ref{ghent7} shows that $H_K(\varepsilon_p, F_p)$ splits.
\end{proof}
As  applications of this theorem, we  used the SageMath software \cite{St} to find all prime Fibonacci numbers $F_p$ ($p> 5$) known until today such that $H_K(\varepsilon_p, F_p)$ and $p$  satisfy conditions 1, 2, 3 and 4 above:
\begin{enumerate}[\indent\rm1.]
	\item $p\in \{13, 433, 25561, 30757, 50833, 130021, 590041, 593689, 1968721\}$.
    \item $p\in \{17, 137, 449, 569, 1049897, 2904353\}$.
    \item $p\in \{11, 131, 359, 431, 9311, 35999, 37511, 81839, 148091, 397379, 1803059\}$.
    \item $p\in \{ 7, 23, 43, 47, 83, 4723, 5387, 201107, 1285607, 1636007\}$.
\end{enumerate}

\end{document}